\documentclass[11pt,reqno]{amsart}

\usepackage{amsmath}
\usepackage{amssymb}
\usepackage{amsthm}
\usepackage{amsbsy}
\usepackage{amsfonts}

\usepackage[bitstream-charter]{mathdesign}
\usepackage[T1]{fontenc}

\newtheorem{thm}{Theorem}
\newtheorem{lem}[thm]{Lemma}
\newtheorem{prop}[thm]{Proposition}

\newtheorem{defn}{Definition}
\newtheorem{rem}[thm]{Remark}
\newtheorem{nota}{Notation}

\usepackage{esint}

\newcommand{\N}{\mathbb{N}}
\newcommand{\R}{\mathbb{R}}
\newcommand{\Z}{\mathbb{Z}}

\newcommand{\cC}{\mathcal{C}}
\newcommand{\cD}{\mathcal{D}}
\newcommand{\cL}{\mathcal{L}}

\newcommand{\m}{\boldsymbol m}

\newcommand{\bj}{\boldsymbol j}

\newcommand{\vep}{\varepsilon}

\newcommand{\intzi}{\int_{0}^{\infty}}
\newcommand{\wL}{L^{3,w}(\R^{3})}

\DeclareMathOperator*{\diam}{diam}
\DeclareMathOperator*{\divv}{div}
\DeclareMathOperator*{\esssup}{ess\,sup}
\DeclareMathOperator*{\loc}{loc}

\newcommand{\set}[1]{\{#1\}}
\newcommand{\Set}[1]{\left\{#1\right\}}
\newcommand{\inn}[1]{\langle#1\rangle}
\newcommand{\norm}[1]{\Vert#1\Vert}

\begin{document}
\baselineskip=18pt

\title[Regularity of the Navier--Stokes equations]{On regularity and singularity for $L^\infty(0,T;L^{3,w}(\mathbb{R}^3))$ solutions to the Navier--Stokes equations}

\author{Hi Jun Choe \& J\"org Wolf \& Minsuk Yang}

\address{H. J. Choe: Department of Mathematics, 
Yonsei University,
Yonseiro 50, Seodaemungu Seoul, Korea}
\email{choe@yonsei.ac.kr}

\address{J. Wolf: Department of Mathematics,
Humboldt University Berlin,
Unter den Linden 6, 10099 Berlin,
Germany}
\email{jwolf@math.hu-berlin.de}

\address{M. Yang: Korea Institute for Advanced Study \\
Hoegiro 85, Dongdaemungu Seoul, Korea}
\email{yangm@kias.re.kr}

\begin{abstract}
We study local regularity properties of a weak solution $u$ to the Cauchy problem of the incompressible Navier--Stokes equations.
We present a new regularity criterion for the weak solution $u$ satisfying the condition $L^\infty(0,T;L^{3,w}(\mathbb{R}^3))$ without any smallness assumption on that scale, where $L^{3,w}(\mathbb{R}^3)$ denotes the standard weak Lebesgue space.
As an application, we conclude that there are at most a finite number of blowup points at any singular time $t$.
The condition that the weak Lebesgue space norm of the veclocity field $u$ is bounded in time is encompassing type I singularity and significantly weaker than the end point case of the so-called Ladyzhenskaya--Prodi--Serrin condition proved by Escauriaza--Sergin--\v{S}ver\'{a}k.
\noindent{\bf keywords:} Navier--Stokes equations, weak Lebesgue space, local regularity
\end{abstract}

\maketitle

\section{Introduction}
\label{S1}

In this paper, we consider the Cauchy problem for the incompressible Navier--Stokes equations
\begin{equation}
\label{E11}
\begin{split}
(\partial_t - \Delta) u+ (u\cdot \nabla) u+ \nabla p = 0 \\
\divv u= 0
\end{split}
\end{equation}
in $Q_T := \R^3\times(0,T)$ and $T>0$ with a smooth and rapidly decaying solenoidal initial vector field $u(x,0)=u_0(x)$ in $\R^3$.
The state variables $u$ and $p$ denote the velocity field of the fluid and its pressure. 
Leray \cite{MR1555394} proved that the Cauchy problem has a unique smooth solution for a short time.
He also proved that there exists at least one global weak solution satisfying an energy inequality.
Hopf \cite{MR0050423} extended the result in the case of bounded domains with a modern concept of weakly differentiable functions.
The weak solution $u$ lies in the space
\begin{equation}
\label{E12}
V^2_\sigma (Q_T) := L^{ \infty}(0,T; L_\sigma^2(\R^3)) \cap L^2(0,T;W_\sigma^{1, 2}(\R^3)),
\end{equation}
but uniqueness and regularity of the weak solution are still open problems.
The exact concept of weak solutions and notations will be given in the next section.

Since there are plenty of important contributions for the regularity question of the Navier--Stokes equations, we briefly describe a few of them closely related to our results.
To guarantee the regularity of weak solutions, one of the most important conditions is the so-called Ladyzhenskaya--Prodi--Serrin \cite{MR0236541,MR0126088,MR0150444} condition, that is,
\begin{equation}
\label{E13a}
u \in L^l(0,T;L^s(\R^3))
\end{equation}
for some $s$ and $l$ satisfying 
\[\frac{3}{s} + \frac{2}{l} = 1, \quad 3 < s \le \infty.\]
Under this condition, the weak solution $u$ to the Cauchy problem \eqref{E11} is unique and smooth.
Later, Escauriaza--Sergin--\v{S}ver\'{a}k \cite{MR2005639} proved that the regularity of a weak solution can also be assured by the marginal case,
\begin{equation}
\label{E13b}
u \in L^\infty(0,T;L^3(\R^3)).
\end{equation}
However, we do not know yet that kinds of higher integrability hold for weak solutions.
By standard embeddings of the solution space \eqref{E12}, any weak solution satisfy the mixed integrability condition with the range of integrability exponents
\[\frac{3}{s} + \frac{2}{l} = \frac{3}{2}, \quad 2 \le s \le 6.\]
There is a considerable gap compared with the Ladyzhenskaya--Prodi--Serrin condition.

To guarantee the local regularity of weak solutions, there are other conditions the so called $\vep$ regularity conditions.
For the Cauchy problem of the Navier--Stokes equations, there is a natural scaling structure 
\begin{align*}
u(x,t) &\to \lambda u(\lambda x, \lambda^2 t), \\
p(x,t) &\to \lambda^2 p(\lambda x, \lambda^2 t).
\end{align*}
Many of the local regularity results have been established under the various smallness assumptions on some scaling invariant quantities. 
We denote by $\Sigma$ the set of possible singular points for the weak solution $u$. 
Utilizing regularity criteria, one can estimate the size of $\Sigma$ by means of some fractal measures and extract some geometric information of $\Sigma$.
In this direction, Scheffer \cite{MR0454426,MR0510154} introduced the concept of suitable weak solutions for the Navier--Stokes equations and then gave partial regularity results.
Caffarelli--Kohn--Nirenberg \cite{MR673830} further strengthened Scheffer's results and gave an improved bound for the Hausdorff dimension of $\Sigma$.
Lin \cite{MR1488514} presented a greatly simplified proof.
Ladyzhenskaya--Seregin \cite{MR1738171} gave more details and considered the case that external forces lie in some Morrey spaces. 
Choe--Lewis \cite{MR1780481} presented an improved estimate of $\Sigma$ in terms of general Hausdorff measures.
Gustafson--Kang--Tsai \cite{MR2308753} unified several known regularity criteria.
For the case \eqref{E13b}, Neustupa \cite{MR1738173} investigated the structure of $\Sigma$ and then Escauriaza--Sergin--\v{S}ver\'{a}k \cite{MR2005639} resolved the regularity qestion.

In this paper, we shall present a new regularity criterion for weak solutions to the Cauchy problem \eqref{E11} satisfying the condition 
\begin{equation}
\label{E14}
u \in L^\infty(0,T;\wL)
\end{equation}
where $\wL$ denotes the weak Lebesgue space.
Because the condition \eqref{E14} is significantly weaker than the condition \eqref{E13b} encompassing type I singularity, the regularity qestion under that condition draws many mathematicians' attention.
However, in the authors knowledge, all results were established under the smallness assumption on that scale \eqref{E14}. 
See, for example, \cite{MR1632780,MR1877269,WZ2014} and the references therein .

We shall use the following notation.

\begin{nota}
We denote the space ball of radius $r$ and center $x$ by $B(x,r) := \set{y \in \R^3 : |y-x|<r}$ and the space-time cylinder at $z = (x,t)$ by 
\[Q(z,r) := B(x,r) \times (t-r^2, t).\]
If the center is at the origin, we simply put $B_r = B(0,r)$ and $Q_r = Q(0,r)$.
\end{nota}

The following theorem is our new regularity criterion.

\begin{thm}
\label{MR1}
For each $M>0$ there exists a positive number $\vep(M) < 1/4$ such that if a weak solution $u\in V_\sigma^2(Q_T)$ to the Cauchy problem \eqref{E11} satisfies the conditions
\begin{equation}
\label{E15}
\esssup_{0 \le t \le T} \norm{u}_{\wL} \le M
\end{equation}
and for some $z_0  = (x_0, t_0) \in Q_T$ and $0 < r \le \sqrt{t_0}$
\begin{equation}
\label{E16}
\frac{1}{r^3} \m \Set{x\in B(x_0,r) : |u(x, t_0)| > \frac{\vep}{r}} \le \vep,
\end{equation}
where $\m(E)$ denotes the Lebesgue measure of the set $E$, then $u$ is bounded in the space-time cylinder $Q(z_0,\vep r)$.
\end{thm}

As an application of this criterion, we are able to estimate the size of possible singular points at a singular time $t$, denoted by 
\[\Sigma(t) = \set{x : (x,t) \in \Sigma}.\]
We know that the Hausdorff dimension of the possible singular time is at most $1/2$.
Many researchers have been investigating the size of $\Sigma(t)$ at the singular time $t$ under various conditions on $u$.
Seregin \cite{MR1829531} obtained a result on estimating $\Sigma(t)$ under the slightly weaker condition than \eqref{E13b}. 
More recently, Wang--Zhang \cite{WZ2014} gave a unifying results on the number of singular points under the Ladyzhenskaya--Prodi--Serrin type conditions.

Utilizing Theorem \ref{MR1}, we can obtain the following theorem which shows that the number of possible singular points at any singular time $t$ is at most finite.

\begin{thm}
\label{MR2}
Suppose $u\in V_\sigma^2(Q_T)$ is a weak solution to the Cauchy problem \eqref{E11} and satisfies
the condition \eqref{E15} for some $M>0$.
Then there exist at most finite number $N(M)$ of singular points at any singular time $t$.
\end{thm}

At each singular time $t$, only a few singular points exist, yet we do not know that blowup points are of type I or not.

\section{Preliminaries}
\label{S2}

Throughout the paper, we shall use the following notation.

\begin{nota}
We denote $A \lesssim B$ if there exists a generic positive constant $C$ such that $|A| \le C|B|$.
We denote the average value of $f$ over the set $E$ by 
\begin{equation}
\inn{f}_E := \fint_E f := \frac{1}{\m(E)} \int_E f
\end{equation}
where $\m(E)$ denotes the Lebesgue measure of the set $E$.
We shall use the same notation $\m$ for the space sets in $\R^3$ and the space-time sets $\R^3 \times (0,T)$ and it will be clearly understood in the contexts.
\end{nota}

We now recall the definition of the weak Lebesgue spaces.
For a measurable function $f$ on $\R^3$, its level set with the height $h$ is denoted by
\begin{equation}
\label{E21}
E(h) = \set{x \in \R^3:|f(x)|>h}.
\end{equation}
The Lebesgue integral can be expressed by the Riemann integral of such level sets.
In particular, for $0<q<\infty$
\begin{equation}
\label{E22}
\int |f(x)|^q dx = \intzi qh^{q-1} \m(E(h)) dh.
\end{equation}

\begin{defn}
The weak Lebesgue space $L^{q,w}(\R^3)$ is the set of all measurable function such that the quantity
\begin{equation}
\label{E23}
\norm{f}_{q,w} := \sup_{h>0} \big[h \m(E(h))^{1/q}\big]
\end{equation}
is finite.
\end{defn}

As the usual convention, two functions are considered the same if they are equal almost everywhere.
In fact, $\norm{f}_{p,w}$ is not a true norm since the triangle inequality fails.
But, it is easy to see that for any $0<r<q$ the following expression
\begin{equation}
\label{E24}
\norm{f} := \sup_{0<\m(E)<\infty} \m(E)^{1/q} \left(\fint_E |f|^r\right)^{1/r}
\end{equation}
is comparable to $\norm{f}_{q,w}$, (see \cite{MR2449250} for example).
Moreover, $\norm{f}$ satisfies the triangle inequality if $1 \le r < q$ and hence it plays the role of true norm for $q>1$.
Furthermore, the weak Lebesgue spaces are Banach spaces and coincide with the Lorentz (Marcinkiewicz) spaces $L^{q,\infty}$.

\begin{rem}
Using \eqref{E24} one can easily see that
\[L^{q,w}(\R^3) \subset \bigcap_{1 \le r < q} \cL_{\loc}^{r,3-3r/q}(\R^3).\]
where $\cL_{\loc}^{r,\lambda}(\R^3)$ denotes the local Morrey space.
\end{rem}

The next remark shows that there is no nonzero harmonic function in $\wL$.
This fact will be used in the proof of our main theorem.

\begin{rem}
\label{R22}
It is easy to see that if $f \in \wL$ is harmonic, then $f=0$.
Indeed, using the mean value property, we have for all $x_0 \in \R^3$ and $R>0$
\begin{align*}
|\nabla f(x_0)| 
&\lesssim \frac{1}{R^4} \int_{B(x_0,R)} |f(x)| dx \\
&\lesssim \frac{1}{R^4} \intzi \m(B(x_0,R) \cap E(h)) dh \\
&\lesssim \frac{1}{R^4} \left(\int_0^H R^3 dh + \int_H^\infty h^{-3} dh\right)
\end{align*}
for all $H>0$.
Taking $H=R^{-1}$ we get $|\nabla f(x_0)| \lesssim R^{-2}$ for all $x_0 \in \R^3$ and $R>0$.
Letting $R \to \infty$, we conclude that $\nabla f=0$ and hence $f$ is constant.
Since $f \in \wL$, it should be identically zero.
\end{rem}

We denote by $L^q(\R^3)$ and $W^{k,q}(\R^3)$ the standard Lebesgue and Sobolev spaces, and we omit these standard definitions.
We denote by $\cD_\sigma(\R^3)$ the set of all solenoidal vector fields $\phi \in C_c^\infty(\R^3)$.
We define $L_\sigma^2(\R^3)$ to be the closure of $\cD_\sigma(\R^3)$ in $L^2(\R^3)$ and $W_\sigma^{1,2}(\R^3)$ to be the closure of $\cD_\sigma(\R^3)$ in $W^{1,2}(\R^3)$.

We now  recall  the concept of local pressure projection (cf. \cite{Wolf2016}). Given a bounded 
$ C^2$-domain $ G \subset \R^{n}$, $ n\in \R^{n}$, we define the operator  
\begin{equation}
E^{\ast}_{G}: W^{-1,\, s}(G) \rightarrow  W^{-1,\, s}(G).
\end{equation}
Appealing to the $ L^p$-theory of the steady Stokes system (cf.  \cite{gal}), 
for any  $ F\in W^{-1,\, s}(G)$ there exists a unique pair $ (v, p)\in W^{1,\, s}_{ 0, \sigma }(G)\times L^s_0(G)$ 
 which solves in the weak sense  the steady Stokes system 
\begin{align*}
-\Delta v + \nabla p = F \quad &\text{ in } \quad  G, \\
\divv v = 0 \quad &\text{ in } \quad  G, \\
v=0 \quad &\text{ on } \quad \partial G.
 \end{align*}
Then we set $E^{\ast}_G(F):=  \nabla p $, where $ \nabla p $ denotes the gradient functional in $ W^{-1,\, s}(G)$  defined by 
\begin{equation}
 \langle \nabla p, \varphi \rangle = \int_{G} p \nabla \cdot \varphi dx, \quad  \varphi \in W^{1,\, s'}_0(G). 
\end{equation}
Here we have denoted by $ L^s_0(G)$ the space of all $ f\in L^s(G)$ with $ \int_{G} fdx=0 $.

\begin{rem}
\label{rem9}
1. The operator $ E^{\ast}_G$  
is bounded from $ W^{-1,\, s}(G)$ into itself with  $ E^{\ast}_G(\nabla p)=\nabla p$ for all $ p\in L^s_{ 0}(G)$. The norm of $ E^{\ast}_G $ depends only on $ s$ and the geometric properties of $ G$, and independent on $ G$,  if $ G$ is a ball or an annulus,  which is due to the scaling properties of the Stokes equation.

2. In case $ F\in L^s(G)$ using the canonical embedding  $ L^s(G) \hookrightarrow  W^{-1,\, s}(G)$ and  the  elliptic regularity 
 we get $ E^{\ast}_G(F)= \nabla p \in L^s(G)$ together with the estimate 
\begin{equation}
\| \nabla p\|_{s, G} \le c \| F\|_{ s, G}, 
\label{E.25a}
\end{equation}
where the  constant in \eqref{E.25a} depends only on $ s$ and $ G$. In case  $ G$  is a ball or an annulus this constant depends only on $ s$ 
 (cf. \cite{gal} for more details). Accordingly the restriction of  $ E^{\ast}_G$ to the Lebesgue space $ L^s(G)$ defines  
a projection in $L^s(G)$. This projection will be denoted still by $ E^{\ast}_G$. 
\end{rem}

Next, we introduce the notion of weak solutions and   local suitable weak solutions.  

\begin{nota}
We denote by $dz$ the space-time Lebesgue measure $dxdt$. 
\end{nota}

\begin{defn}
We say that $u$ is a Leray--Hopf weak solution to \eqref{E11} if the velocity field $u$ lies in the space $V^2_\sigma (Q_T) = L^{ \infty}(0,T; L_\sigma^2(\R^3)) \cap L^2(0,T;W_\sigma^{1, 2}(\R^3))$, there exists a distribution $p$ such that $(u, p)$ solves the Navier--Stokes equations in the sense of distributions, and $u$ satisfies the energy inequality for almost all $s \in (0,T)$
\[\int_{\R^3} |u(s)|^2 dx + 2 \int_0^s \int_{\R^3} |\nabla u|^2 dz 
\le \int_{\R^3} |u(0)|^2 dx.\]
We say that $u $ is a local suitable weak solution to \eqref{E11} if 
for every ball $ B \subset \R^{3}$
the following local energy inequality  the following local energy inequality holds for almost all $ s \in (0,T)$ and  for all non negative $\phi \in C^{\infty}_{\rm c}(B\times (0,T))$, 
\begin{equation}
\begin{split}
\label{E.25}
&\int |v(s)|^2 \phi (s) dx + 2 \int_{0}^{s} \int| \nabla v|^2  \phi dz \\
&\le \int_{0}^{s} \int |v|^2 (\partial_t + \Delta) \phi dz
+ \int_{0}^{s} \int |v|^2 (v - \nabla p_h) \cdot \nabla \phi dz \\
&\quad + 2 \int_{0}^{s} \int (v \otimes v - v \otimes \nabla p_h  : \nabla ^2 p_h) \phi dz
+ 2 \int_{0}^{s}\int (p_{ 1, B}+ p_{ 2, B}) v \cdot \nabla \phi dz,
\end{split}
\end{equation}
where $ v= u+ \nabla p_{ h,B}$, and 
\begin{align*}
\nabla p_{h,B} &= - E^{\ast}_{ B} (u), \\
\nabla p_{1,B} &= - E^{\ast}_{ B} (\nabla \cdot (u \otimes u)), \\
\nabla p_{2,B} &= E^{\ast}_{ B}(\Delta u).
\end{align*}
\end{defn}

\begin{rem}
If a weak solution $u$ is in $L^\infty(0,T;\wL)$, then $u$ lies in $L^4(Q_T)$ by an interpolation.
Thus, the function $|u|^2|\nabla u|$ is integrable on $Q_T$, which justifies the integration by parts and one can show that $u$ becomes a local suitable weak solution, too.
\end{rem}

Using a standard iteration method one can observe that boundedness of a certain scaling invariant quantity essentially implies the boundedness of many of other scaling invariant quantities.
The following form of the Caccioppoli-type inequality is convenient in that purpose.

\begin{lem}[Lemma\,2.6 in \cite{ChaeWolf2016a}]
\label{L22}
If $u$ is a suitable weak solution to \eqref{E11}, then for all $Q(z_0,r) \subset Q_T$
\begin{equation}
\label{E26}
\begin{split}
&r^{-1} \left(\int_{Q(z_0,r/2)} |u|^{10/3} dz\right)^{3/5}
+ r^{-1} \int_{Q(z_0,r/2)} |\nabla u|^2 dz \\
&\lesssim \left(r^{-5} \int_{t_0-r^2}^{t_0} \left(\int_{B(x_0,r)} |u|^2 dx \right)^3 dt\right)^{1/3} 
+ r^{-5} \int_{t_0-r^2}^{t_0} \left(\int_{B(x_0,r)} |u|^2 dx \right)^3 dt
\end{split}
\end{equation}
where the implied constant is absolute.
\end{lem}

We end this section by giving the following version of the local regularity criterion.
We include its proof at the end of this paper, Appendix \ref{SA1}.

\begin{lem}[\cite{Wolf2015c}]
\label{L22}
There exists an absolute positive number $\zeta$ such that if a local suitable weak solution $ u\in 
V^2_\sigma (Q(z_0,\rho))$ to the Navier--Stokes equations satisfies the condition
\begin{equation}
\label{E27}
\rho^{-2} \int_{Q(z_0,\rho)} |u|^3 dz \le \zeta^3
\end{equation}
for some $z_0  = (x_0, t_0)\in Q_T$ and $ 0 < \rho \le \sqrt{t_0}$, then 
\[u \in L^\infty(Q(z_0,\rho/2)), 
\]
and  the following estimate holds true 
\begin{equation}
\| u\|_{ L^\infty(Q(z_0,\rho/2))} \le  C\left(\fint_{Q(z_0,\rho)} | u|^3 dz\right)^{1/3} + 
C \esssup_{t\in (t_0-\rho ^2,t_0)} \left(\fint_{B(x_0,\rho)} | u(t)|^2 dx\right)^{1/2}
\label{E.28}
\end{equation}
where $C$ is an absolute positive constant. 
\end{lem}

\section{Proof of Theorem \ref{MR1}}
\label{S3}

Due to Lemma \ref{L22}, it suffices to show that the following lemma holds true.

\begin{lem}
\label{L31}
For each $M>0$ there exists a positive number $\vep(M) < 1/4$ such that if a weak solution $u\in V_\sigma^2(Q_T)$ to the Navier--Stokes equations satisfies the condition 
\begin{equation}
\label{E31}
\esssup_{0 \le t \le T} \norm{u}_{\wL} \le M
\end{equation}
and for some $z_0  = (x_0, t_0) \in Q_T$ and $0 < r \le \sqrt{t_0}$
\begin{equation}
\label{E32}
r^{-3} \m \set{x\in B_r(x_0) : |u(x, t_0)| > r^{-1} \vep} \le \vep,
\end{equation}
then there exists $\rho \in [2\vep r, \sqrt{t_0}]$ such that
\begin{equation}
\label{E33}
\rho^{-2} \int_{Q(z_0,\rho)} |u|^3 dz \le \zeta^3
\end{equation}
where $\zeta$ is the same number in Lemma \ref{L22}.
\end{lem}

We divide the proof of Lemma \ref{L31} into several steps.

\begin{itemize}
\item[\textbf{Step 1)}]
We first observe that the condition \eqref{E31} yields 
\[u \in C([0,T]; L^2(\R^{3})).\]
Indeed, \eqref{E31} implies that for almost all $0 \le t \le T$ and all $h>0$
\begin{equation}
\label{E34}
h^3 \m(E_t(h)) \le M^3
\end{equation}
where $E_t(h)$ denotes the level set 
\[E_t(h) := \set{x:|u(x,t)|>h}.\]
By the Chebyshev inequality we also have 
\begin{equation}
\label{E35}
h^6 \m(E_t(h)) \le \norm{u(t)}_{L^6}^6.
\end{equation}
Using the two estimates \eqref{E34} and \eqref{E35}, we obtain that for any $H>0$
\begin{align*}
\int_{\R^3} |u(x,t)|^4 dx 
&= 4 \intzi h^3 \m(E_t(h)) dh \\
&\lesssim \int_0^H M^3 dh + \int_H^\infty h^{-3} \norm{u(t)}_{L^6}^6 dh \\
&\lesssim M^3H + \norm{u(t)}_{L^6}^6H^{-2}.
\end{align*}
Taking $H = M^{-1}\norm{u(t)}_{L^6}^2$ we get 
\[\int_{\R^3} |u(x,t)|^4 dx \lesssim M^2 \norm{u(t)}_{L^6}^2.\]
Hence $u \in L^4(Q_T)$ and so $|u|^2|\nabla u| \in L^1(Q_T)$.
This justifies the required integration by parts to be a local suitable weak solution and also implies the global energy equality so that $u$ is in $C([0,T];L^2(\R^3))$.

\item[\textbf{Step 2)}]
We next claim that the condition \eqref{E31} also yields that for all $Q(z_0,r) \subset Q_T$
\begin{equation}
\label{E36}
r^{-1} \left(\int_{Q(z_0,r/2)} |u|^{10/3} dz\right)^{3/5}
+ r^{-1} \int_{Q(z_0,r/2)} |\nabla u|^2 dz  
\lesssim M^2+M^6.
\end{equation} 
Due to the Caccioppoli--type inequality \eqref{E26}, it suffices to estimate
\[\int_{B(x_0,r)} |u(x,t)|^2 dx.\] 
Using the estimate \eqref{E34}, we obtain that for almost all $0 \le t \le T$ and all $h>0$
\begin{align*}
\int_{B(x_0,r)} |u(x,t)|^2 dx 
&= 2 \intzi h \m [B(x_0,r) \cap E_t(h)] dh \\
&\lesssim \int_0^H h r^3 dh + \int_H^\infty h^{-2} M^3 dh \\
&\lesssim r^3 H^2 + M^3 H^{-1}.
\end{align*}
Taking $H = Mr^{-1}$ we get for almost all $0 \le t \le T$
\begin{equation}
\label{E37}
\int_{B(x_0,r)} |u(x,t)|^2 dx \lesssim M^2 r.
\end{equation}
Putting this bound into the right side of the inequality \eqref{E26}, we get the estimate \eqref{E36}.

\item[\textbf{Step 3)}]
We now prove Lemma \ref{L31} by using an indirect argument.
Assume the assertion of the lemma is not true, that is, there exist a positive number $M$, sequences $\vep_k \in (0, 1/4)$, $T_k \in (0,\infty)$, $z_k = (x_k,t_k) \in Q_{T_k}$, $r_k \in (0,\sqrt{t_k}]$, and a sequence of weak solutions $u_k \in V^2_\sigma (Q_{ T_k})$ such that $\vep_k \to 0$ as $k \to \infty$ and for all $k \in \N$
\begin{equation}
\begin{split}
\label{E38}
\esssup_{0 \le t \le T_k} \norm{u_k}_{\wL} &\le M, \\
r_k^{-3} \m \set{x \in B(x_k,r_k) : |u_k(x,t_k)| > r_k^{ -1} \vep_k} &\le \vep_k, 
\end{split}
\end{equation}
and for all $\rho \in (2\vep_k r_k, \sqrt{t_k}]$
\begin{equation}
\label{E39}
\rho^{-2} \int_{Q(z_k,\rho)} |u_k|^3 dz > \zeta^3.
\end{equation}
We define for $(y,s) \in \R^3 \times (-1,0)$
\begin{align*}
U_k(y,s) &= r_k u_k(x_k+ r_k y, t_k+ r_k^2 s), \\
P_k(y,s) &= r^2_k p_k(x_k+ r_k y, t_k+ r_k^2 s).
\end{align*}
Then $(U_k, P_k)$ is a weak solution to the Navier--Stokes equations in $\R^{3}\times (-1,0)$ and satisfies
\[\esssup_{-1 \le s \le 0} \norm{U_k}_{\wL} \le M.\]
Thanks to \eqref{E36} and \eqref{E37}, we have for all $ k\in \N$, $z_0 = (x_0,0)$ and $0 < \rho \le 1$
\begin{equation}
\label{E310}
\rho^{-1} \left(\int_{Q(z_0,\rho/2)} |U_k|^{10/3} dz\right)^{3/5}
+ \rho^{-1} \int_{Q(z_0,\rho/2)} |\nabla U_k|^2 dz \lesssim M^2+M^6
\end{equation}
and 
\begin{equation}
\label{E311}
\rho^{-1} \sup_{-\rho^2 \le s \le 0} \int_{B(x_0,\rho)} |U_k(s)|^2 dx \lesssim M^2.
\end{equation}
Furthermore, from \eqref{E38} and \eqref{E39}, we also have for all $ k\in \N$
\begin{equation}
\label{E312}
\m \set{x \in B_1 : |U_k(x,0)| > \vep_k} \le \vep_k
\end{equation}
and for all $\rho \in [\vep_k,1]$
\begin{equation}
\label{E313}
\rho^{-2} \int_{Q(0,\rho)} |U_k|^3 dz  > \zeta^3.
\end{equation}
Using a standard reflexivity argument along with Cantor's diagonalization principle and passing to a subsequence from \eqref{E310} we eventually get $U \in L^{10/3}(-1,0; L^{10/3}_{\loc}(\R^{3}))$ with $ \nabla U \in L^{2}(-1,0; L^{2}_{\loc}(\R^{3}))$ and $H\in L^{5/3}(-1,0; L^{5/3}_{\loc}(\R^{3}))$ such that for every $ 0<R<\infty$ 
\begin{equation}
\label{E314}
\begin{split}
U_k \rightarrow U \quad &\text{{\it weakly in}} \quad L^{10/3}(B_R \times (-1,0)) \\
\nabla U_k \rightarrow \nabla U \quad &\text{{\it weakly in}} \quad L^{2}(B_R\times (-1,0)) \\
U_k \otimes U_k \rightarrow H  \quad &\text{{\it weakly in}}\quad  L^{5/3}(B_R\times (-1,0))
\end{split}
\end{equation}
as $k \rightarrow \infty$.
Hence, $ U$ appears to be a distributional solution to 
\begin{equation}
\label{E326}
\partial_t U - \Delta  U + \nabla \cdot H = -\nabla P \quad  \text{ in}\quad  \R^{3}\times (-1,0). 
\end{equation}   
According to the weakly lower semi-continuity of the norm we get from \eqref{E310} and \eqref{E311}  along with  \eqref{E314} 
for all $0 < \rho \le 1$
\begin{equation}
\label{E327}
\rho^{-1} \left(\int_{Q(z_0,\rho)} |U|^{10/3} dyds\right)^{3/5} 
+ \rho^{-1} \int_{Q(z_0,\rho)} | \nabla U|^{2} dy ds \lesssim M^2 + M^6
\end{equation}
and 
\begin{equation}
\label{E328}
\rho^{-1} \sup_{-\rho^2 \le s \le 0} \int_{B(x_0,\rho)} |U(s)|^2 dy \lesssim M^2.
\end{equation}

\item[\textbf{Step 4)}]
Let $s_0 \in [-1,0]$.  Since $u \in C_{ w}^{\ast}([0,T]; \wL)$, we have 
$u(\cdot , t_k + r_k^2 s_0) \in \wL$ and 
\[\norm{u(\cdot , t_k + r_k^2 s_0)}_{ \wL} = \norm{U_k(\cdot,s_0)}_{ \wL}.\]
This shows that $\{U_k(\cdot,s_0)\}$ is a bounded sequence in $\wL$. 
On the other hand, the predual of $ \wL$ is the Lorentz space $ L^{3/2,1}(\R^{3})$. 
By means of the Banach--Alaoglu theorem we get a subsequence 
$ \{ U_{ k_j}(\cdot ,  s_0)\}$  and a function $ \eta \in \wL$ such that 
\begin{equation}
 U_{ k_j}(\cdot ,  s_0) \rightarrow \eta \quad  \text{{\it weakly$^{\ast}$ in}}\quad \wL\quad  \text{{\it as}}\quad  j \rightarrow +\infty. 
\label{E329}
\end{equation}
Thus, from \eqref{E326} we infer that for all $ \varphi \in C_c^{\infty}(\R^{3}\times (-1,0)) $ with $ \nabla \cdot \varphi =0$
\[\int_{-1}^{s_0} \int_{ \R^{3}} - U \cdot \partial_t \varphi  + \nabla U: \nabla \varphi  - H : \nabla \varphi dz
 = - \int_{ \R^{3}}  \eta \cdot \varphi (s_0) dx.
\] 
In case $ s_0$ is a Lebesgue point of $U$ with respect to time, we argue that for all $ \psi \in C^{\infty}_{\rm c, \sigma }(\R^{3})$
\[
\int_{\R^{3}} (U(s_0) - \eta ) \cdot \psi dx=0 
\]
which shows that $ U(s_0) - \eta$ is a gradient field. Together with $ \nabla \cdot (U(s_0) - \eta)$ in the sense of distributions 
we see that $ U(s_0) - \eta$ is harmonic. Recalling that $ U(s_0) - \eta \in \wL$, it follows that $ \eta = U(s_0)$ by Remark \ref{R22} in Section \ref{S2}. 
Consequently, \eqref{E329} yields  
\[U_k(\cdot ,  s_0) \rightarrow U(s_0) \quad  \text{{\it weakly$^{\ast}$ in}} \quad \wL \quad  \text{{\it as}} \quad k \rightarrow +\infty.\]
Furthermore, we get 
\[\norm{U(s_0)}_{\wL} \le \norm{u}_{L^\infty(0,T;\wL)}.\]
In particular, $U \in L^\infty(0,T;\wL)$. 

\item[\textbf{Step 5)}]
Next, let $ s_0\in [-1,0]$. Then we may choose $ s_m \in (0,T)$ in the set of Lebesgue points such that $ s_m \rightarrow s_0$ as $ m \rightarrow +\infty$.  Then as above we get a subsequence $ \{ s_{ m_j}\}$ and $ \eta \in \wL$ such that 
\[U(\cdot ,  s_{ m_j}) \rightarrow \eta \quad  \text{{\it weakly$^{\ast}$ in}}\quad \wL\quad  \text{{\it as}}\quad  j \rightarrow +\infty.\]
In addition, we easily verify that the following identity holds for every 
$ \varphi \in C^{\infty}(\R^{3}\times (-1,0)) $ with $\divv \varphi = 0$
\begin{equation}
 \int_{-1}^{s_0} \int_{ \R^{3}} - U \cdot  \partial_t \varphi  + \nabla U: \nabla \varphi  - H : \nabla \varphi dz
 = - \int_{ \R^{3}}  \eta \cdot \varphi (s_0) dx.  
\label{E331}
\end{equation}
Arguing as above, we see that this limit is unique, and will be denoted by $ U(s_0)$. Note that  \eqref{E331} holds true for $ \eta = U(s_0)$.  
We now repeat the same argument as above to prove that for all $ s_0 \in [-1,0]$
\begin{align}
\label{E332}
U_k(\cdot,s_{ 0}) &\rightarrow U(s_0) \quad  \text{{\it weakly$^{\ast}$ in}}\quad \wL\quad  \text{{\it as}}\quad  k\rightarrow +\infty \\
\label{E333}
U(\cdot,s) &\rightarrow U(s_0) \quad  \text{{\it weakly$^{\ast}$ in}}\quad \wL\quad  \text{{\it as}}\quad  s\rightarrow s_0.  
\end{align}
This leads to $ U\in C_{ w}^{\ast}([-1,0]; \wL)$. 

\item[\textbf{Step 6)}]
We shall verify the strong convergence of $ U_k$ in $ L^2(B_R\times (-1,0))$. For this purpose,  we define the local pressure 
introduced  in \cite{Wolf2016}, 
\begin{align*}
\nabla P_{h,k,R} &= -E^{\ast}_{B_R} (U_k), \\
\nabla P_{1,k,R} &= -E^{\ast}_{B_R} (\nabla \cdot U_k \otimes U_k), \\
\nabla P_{2,k,R} &= E^{\ast}_{B_R} (\Delta U_k)
\end{align*}
and 
\begin{align*}
\nabla P_{h,R} &= -E^{\ast}_{B_R} (U), \\
\nabla P_{1,R} &= -E^{\ast}_{B_R} (\nabla \cdot H), \\
\nabla P_{2,R} &= E^{\ast}_{B_R} (\Delta U )
\end{align*}
(For the definition of $E^{\ast}_{B_R}$ see Appendix \ref{SA2} of this paper).

\item[\textbf{Step 7)}]
We set $ V_k = U_k + \nabla P_{ h,k,R}$,  and $ V = U+ \nabla P_{ h,R}$. Then $ V_k$ solves 
\[\partial_t V_k - \Delta U_k + \nabla \cdot (U_k \otimes U_k) = -\nabla (P_{ 1,k,R}+ P_{ 2,k,R})
\quad  \text{ in}\quad  B_R\times (-1,0),\]
while $ V$ solves 
\[\partial_t V - \Delta U + \nabla \cdot H = -\nabla (P_{ 1,R}+ P_{ 2,R})\quad  
\text{ in}\quad  B_R \times (-1,0).\]
By using a standard compactness argument due to Lions-Aubin we see that 
\[V_k \rightarrow V  \quad  \text{{\it in}}\quad  L^2(B_R\times (-1,0))\quad  \text{{\it as}}\quad  k \rightarrow +\infty.\]
By passing to a subsequence we may also assume that 
\[V_k \rightarrow V  \quad  \text{{\it a.\,e. in \, $ B_R\times (-1,0)$}}\quad   \text{{\it as}}\quad  k \rightarrow +\infty.\]
Arguing as in \cite{Wolf2015c}, by the aid of \eqref{E333}, and noting  that $ P_{ h, k, R}$ is harmonic, we also find that 
\[\nabla P_{h,k,R } \rightarrow \nabla P_{ h, R}  \quad  \text{{\it a.\,e. in \, $ B_R\times (-1,0)$}}\quad   \text{{\it as}}\quad  k \rightarrow +\infty.\]
This leads to the a.\,e. convergence of $ U_k$ which allows to apply Lebesgue's dominated convergence theorem. 
Accordingly, 
\[U_k \rightarrow U  \quad  \text{{\it in}} \quad L^3(B_R\times (-1,0)) \quad  \text{{\it as}}\quad  k \rightarrow +\infty.\]
This also shows that $ H= U \otimes U$ and therefore $ U$ solves the Navier-Stokes equations. 

\item[\textbf{Step 8)}]
In \eqref{E313} letting $ k \rightarrow +\infty$, we obtain for every $ 0<\rho \le 1$
\begin{equation}
\label{E334}
\rho^{-2} \int_{Q(0,\rho)} |U|^3 dz \ge \zeta^3.
\end{equation} 
It remains to carry out the passage to the limit $ k \rightarrow +\infty$ in \eqref{E312}. Without loss of generality we may assume $\vep_k \le 2^{ -k}$.
Let
\[A := \bigcap_{m=1}^\infty \bigcup_{k=m}^\infty \set{x \in B_1 : |U_k(x,0)| > \vep_k}.\] 
Then according to \eqref{E312} we have
\[\sum_{k=1}^\infty \m \set{x \in B_1 : |U_k(x,0)| > \vep_k} \le \sum_{k=1}^\infty \vep_k < \infty.\]
Hence the Borel--Cantelli lemma yields $\m(A)=0$.
In other words, for each $x \in B_1 \setminus A$, there exists $m \in \N$ such that for all $k \ge m$
\[|U_k(x,0)| \le \vep_k.\]
Accordingly, $U_k(x,0) \to 0$ for almost all $x \in B_1$.
In view of \eqref{E332} we conclude that 
\begin{equation}
\label{E335}
U(0) = 0 \quad \text{ on } \quad  B_1. 
\end{equation} 

\item[\textbf{Step 9)}]
Next, we set $\rho_k = 2^{-k}$ and define for $(x,t) \in \R^{3}\times (-1,0)$
\begin{align*}
{\tilde U}_k(x,t) = \rho_k U(\rho_k x, \rho_k^2 t), \\
\quad {\tilde P}_k(x,t) = \rho^2_k P(\rho_k x,\rho_k^2 t).   
\end{align*}
Again $({\tilde U}_k, {\tilde P}_k)$ is a solution to the Navier-Stokes equation in $\R^{3}\times (-1,0)$.  
Observing \eqref{E327} and \eqref{E328}, we find for all $z_0=(x_0,0)$ 
\begin{equation}
\label{E336}
\left(\int_{Q(z_0,1)} |{\tilde U}_k|^{10/3} dz\right)^{3/5} 
+ \int_{Q(z_0,1)} |\nabla {\tilde U}_k|^{2} dz \lesssim M^2+M^6,
\end{equation}
and 
\begin{equation}
\label{E337}
\sup_{-1 \le t \le 0} \int_{B(x_0,1)} |{\tilde U}_k(t)|^2 dx \lesssim M^2.
\end{equation}
On the other hand, \eqref{E335} and \eqref{E334} yield ${\tilde U} _k(0) = 0$ on $B_{2^k}$ and 
\begin{equation}
\label{E338}
16 \int_{Q(0,1/4)} |{\tilde U}_k|^3 dz > \zeta^3. 
\end{equation}
Arguing as in \textbf{Step 3}, we get a solution 
\[{\tilde U} \in C_{w}^{\ast}([-1,0]; \wL) \cap L^2(-1,0; W^{1,2}_{\loc}(\R^{3}))\]
to the Navier-Stokes equations. 
Furthermore, \eqref{E336}, \eqref{E337}, and \eqref{E338}  yield for all $z_0=(x_0,0)$ 
\begin{equation}
\label{E339}
\left(\int_{Q(z_0,1)} |{\tilde U}|^{10/3} dz\right)^{3/5} 
+ \int_{Q(z_0,1)} |\nabla {\tilde U}|^{2} dz \le C_0,
\end{equation}
\[\sup_{-1 \le t \le 0} \int_{B(x_0,1)} |{\tilde U}(t)|^2 dx \le C_0,\]
${\tilde U} (0) = 0$ in $\R^{3}$, and 
\begin{equation}
\label{E340}
16 \int_{Q(0,1/4)} |{\tilde U}|^3 dyds  \ge  \zeta^3. 
\end{equation}

\item[\textbf{Step 10)}]
By the Fubini theorem, we have
\begin{align*}
&\qquad \m \set{(x,t) \in \R^{3} \times (-1,0) : |{\tilde U}(x,t)| \ge 2^{-5} \zeta} \\
&\qquad = \int_{-1}^{0} \m \set{x \in \R^{3} : |{\tilde U}(x,t)| \ge 2^{-5} \zeta} dt \\
&\qquad \le 2^{15} \zeta ^{ -3} \norm{{\tilde U}}_{L^\infty(-1,0;\wL)}^3 
\\
&\qquad < \infty.
\end{align*}
Hence, for each $\eta>0$ there is a radius $R=R(\eta)>0$ such that 
\[\m \set{(x,t) \in (\R^{3} \setminus B_R) \times (-1,0) : |{\tilde U}(x,t)| \ge 2^{-5} \zeta} \le \eta.\]
Choose 
\begin{equation}
\label{E341}
\eta = 2^{-10} C_0^{-15} \zeta^{30},
\end{equation}
where $C_0$ is the constant in \eqref{E339}.
Then for any $x_0 \in \R^3 \setminus B(0,R+1)$, we obtain, by H\"older's inequality, \eqref{E339}, and \eqref{E341}, that 
\begin{align*}
\int_{Q(z_0,1)} |{\tilde U}|^3 dz 
&\le (2^{-5} \zeta)^3 \m(Q(z_0,1)) +
\int_{Q(z_0,1) \cap \set{|{\tilde U}| \ge 2^{-5} \zeta}} |{\tilde U}|^3  dz 
\\
&\le \frac{\zeta^3}{2} + 
\eta^{1/10} \left(\int_{Q(z_0,1)} |{\tilde U}|^{10/3}  dz\right)^{9/10} 
\\
&\le \frac{\zeta^3}{2} + 2^{-1} C_0^{-3/2} \zeta^{3} C_0^{3/2} 
\\
&\le \zeta^3.   
\end{align*}
Thus, appealing to Lemma\,\ref{L22},   and making use of \eqref{E339}, we get for all $x_0 \in \R^{3} \setminus B(0,R+1)$ 
\begin{equation}
\norm{{\tilde U}}_{L^\infty(Q(z_0,1/2))} \le C\norm{{\tilde U}}_{L^3(Q(z_0,1))} 
+ C \norm{{\tilde U}}_{L^\infty(-1,0;L^2(B(x_0, 1)))}
\le C (\zeta +  C_0).
\end{equation}
This shows that ${\tilde U}$ is bounded in $\R^{3} \setminus B(0,R+1) \times (-1/4,0)$.  
 
\item[\textbf{Step 11)}]
Using a standard bootstrapping argument, we obtain the higher regularity 
\[\nabla {\tilde U} \in L^{\infty} \left(\R^{3} \setminus B(0,R+2) \times (-1/4,0)\right).\] 
Taking the curl operator to the Navier--Stokes equations, we see that ${\tilde \Omega} := \nabla \times {\tilde U}$ solves the heat equation
\[\partial_t {\tilde \Omega} - \Delta {\tilde \Omega } = {\tilde \Omega } \cdot \nabla {\tilde U} - {\tilde U} \cdot \nabla {\tilde \Omega}\]
in $\R^{3}\times (-1,0)$.
Hence 
\[\left|\partial_t {\tilde \Omega} - \Delta {\tilde \Omega}\right| 
\le \norm{\nabla {\tilde U}}_{\infty} |{\tilde \Omega}| 
+ \norm{{\tilde U}}_{\infty} |\nabla {{\tilde \Omega}}|\]
in $\R^{3} \setminus B(0,R+2) \times (-1/16,0)$.
Verifying that $ {\tilde \Omega } (0)=0$,  we are in a position to apply the backward uniqueness of \cite{MR2005639} to conclude that ${\tilde \Omega } \equiv 0$ in $\R^{3} \setminus B(0,R+2) \times 
(-1/16,0)$.
By the spatial analyticity of ${\tilde U}$ we get the spatial analyticity of ${\tilde \Omega } $ which shows that 
$ {\tilde \Omega } \equiv 0  $ in $   \R^{3}  \times  (-1/16,0)$. 
Recalling that $\divv {\tilde U} = 0$ it follows that $ {\tilde U} $ is harmonic in $  \R^{3}\times  (-1/16,0)$, and thus $ {\tilde U} $ must be identically zero in 
$ \R^{3}\times   (-1/16,0)$. 
However this contradicts to \eqref{E340}. 
Therefore the assertion of Lemma \ref{L31} must be true.
\end{itemize}

This completes the proof of Lemma \ref{L31}.
By combining Lemma \ref{L22} we obtain Theorem \ref{MR1}.

\section{Proof of Theorem \ref{MR2}}
\label{S4}

We divide the proof of Theorem \ref{MR2} into a few steps.

\begin{itemize}
\item[\textbf{Step 1)}]
Let $C(x_0,r)$ denote the closed cube of a side-length $r$ and the center $x_0$.
We may replace the condition \eqref{E16} in Theorem \ref{MR1} by using cubes, that is,
\begin{equation}
\label{E44}
r^{-3} \m \set{x\in C(x_0,r) : |u(x, t_0)| > r^{-1} \vep} \le \vep.
\end{equation}
Then the conclusion also be changed with $u\in L^\infty(\widetilde{Q}(z_0,\vep r))$ where 
\[\widetilde{Q}(z_0,\vep r) := C(x_0,r) \times (t_0-r^2,t_0).\] 
In fact, $\vep$ should be changed by a multiplication of some constant which depends only on the volume ratio of the ball of a radius $r$ and the cube of a side-length $r$.
For convenience we just use the same letter $\vep$.

\item[\textbf{Step 2)}]
We shall proceed with an algorithm based on a dyadic decomposition argument.
We say that two cubes $E$ and $E'$ meet if $E \cap E'$ has nonempty interior.
Let $C = [0,1]^3$ denote the unit cube in $\R^3$.
We define for $k=0,1,2,\dots$ the following covers 
\[\cC_k := \set{2^{-k} (\vep \bj + C) : \bj \in \Z^3},\]
which has finite overlapping property.
Indeed, each fixed cube in $\cC_k$ can meet $\vep^{-3}$ number of cubes in $\cC_k$.
We pick a sub-family
\begin{equation}
\label{E45}
F_0 := \set{E \in \cC_0 : \m \set{x \in E : |u(x,t_0)|>\vep} > \vep}.
\end{equation}
If $F_0$ has no element, then we have $\m \set{x \in E : |u(x,t_0)|>\vep} \le \vep$ for all $E \in \cC_0$.
Hence we conclude that there is no singularity at all at the moment $t_0$ due to Theorem \ref{MR1}.

Next, we claim that $F_0$ has at most a finite number of members, which is bounded by a number depending only on $M$ and $\vep$.
Suppose that $E_1, E_2, \dots, E_N \in F_0$ don't meet each other.
Then for $j=1,2,\dots,N$
\[\vep < \m \set{x \in E_j : |u(x,t_0)|>\vep}.\]
Summing both sides for $j=1,2,\dots,N$ yields
\[N\vep < \m \set{x \in \bigcup_{j=1}^N E_j : |u(x,t_0)| > \vep} \le \vep^{-3} M^3.\]
The last inequality follows from the fact $\norm{u(t_0)}_{\wL} \le M$.
This implies that the number of maximal disjoint cubes in $F_0$ is finite and hence $F_0$ has at most finite members.
If we denote by $N_0^d$ the number of maximal disjoint cubes in $F_0$, then we should have 
\[N_0^d \le \vep^{-4} M^3.\]
Let $N_0$ denote the number of cubes in $F_0$.
Then, from the finite overlapping property of $\cC_0$, we have 
\[\vep^3 N_0 \le N_0^d \le N_0.\]
Hence 
\begin{equation}
N_0 \le \vep^{-3} N_0^d \le \vep^{-3} (\vep^{-4} M^3) = \vep^{-7} M^3.
\end{equation}
We define $G_0$ to be the union of $F_0$ and the cubes $E \in \cC_0$ which meet some element of $F_0$.
Theorem \ref{MR1} implies that if $(x,t_0) \notin \bigcup_{E \in G_0} E$, then $(x,t_0)$ is a regular piont, that is, possible singularities can only occur in some element of $G_0$.

\item[\textbf{Step 3)}]
We now inductively construct two families of cubes $\set{F_k}$ and $\set{G_k}$.
For $k\ge1$ we define $F_k$ to be the family of cubes $E \in \cC_k$ satisfying 
$E \subset E'$ for some $E' \in G_{k-1}$ and 
\begin{equation}
\label{E46}
\m \set{E : |u(x,t_0)|>2^k\vep} > 2^{-3k}\vep.
\end{equation}
Let $N_k$ denote the number of cubes in $F_k$ and let $N_k^d$ denote the number of maximal disjoint cubes in $F_k$.
By the same reasoning $N_k$ and $N_k^d$ are finite numbers and have the same bounds.
Indeed, since each fixed cube in $F_k$ can meet at most $\vep^{-3}$ number of cubes in $F_k$, we have 
\[\vep^3 N_k \le N_k^d \le N_k.\] 
By the same way in the previous step, we obtain
\[N_k^d 2^{-3k} \vep \le \m \set{x \in \R^3 : |u(x,t_0)| > 2^k\vep} \le (2^k\vep)^{-3} M^3.\]
Therefore,
\begin{equation}
\label{E47}
N_k \le \vep^{-3} N_k^d \le \vep^{-3} (\vep^{-4} M^3) = \vep^{-7} M^3.
\end{equation}
We define $G_k$ to be the union of $F_k$ and the cubes $E \in \cC_k$ which meet some element of $F_k$.
Theorem \ref{MR1} implies that if $(x,t_0) \notin \bigcup_{E \in G_k} E$, then $(x,t_0)$ is a regular piont, that is,  the possible singularities can only occur in the elements of $G_k$.

\item[\textbf{Step 4)}]
Finally, we construct nested sequences $\set{E_k}$ of closed cubes satisfying $E_k \in G_k$.
Fix an element $E_k$ in $G_k$.
If $E \in G_{k+1}$, then $E \subset E_k$ or $E$ does not meet $E_k$ by the dyadic construction. 
If there is no $E \in G_{k+1}$ which meet $E_k$, then each interior point of $E_k$ is a regular point.
In this case, we stop to choose next elements.
Otherwise, there is an element $E \in G_{k+1}$ such that $E \subset E_k$.
Then we pick $E$ and name it as $E_{k+1}$.
The cardinality of each set $G_k$ is bounded by $\vep^{-7} M^3 + \vep^{-3}$ from \eqref{E47} and the finite overlapping property. 
The number of such choices is also always bounded by $\vep^{-7} M^3 + \vep^{-3}$.
After the construction, we only have at most $\vep^{-7} M^3 + \vep^{-3}$ number of sequences $\set{E_k}$.
If the sequence $\set{E_k}$ is finite, then each interior point of $E_k$ is regular.
If the sequence $\set{E_k}$ is infinite, then 
\[E_{k+1} \subset E_k, \quad \diam E_{k+1} \le \frac{1}{2} \diam E_k.\]
where $\diam E$ denote the diameter of the set $E$.
Since $\diam E_k$ goes to $0$ as $k \to \infty$, 
\[\bigcap_{k=1}^\infty E_k = (x,t_0)\]
for some $x \in \R^3$.
This point might be a singular point.
Therefore, the number of such possible singularities is at most 
\begin{equation}
\label{E48}
N(M) := \vep^{-7} M^3 + \vep^{-3}
\end{equation}
at the time $t_0$.
We note that $\vep$ actually depends only on $M$, and hence the number of possible singularities is bounded by the uniform norm of the weak Lebesgue space.
\end{itemize}
This completes the proof of Theorem \ref{MR2}.

\section*{Acknowledgement}

H. J. Choe has been supported by the National Reserch Foundation of Korea (NRF) grant, funded by the Korea government(MSIP) (No. 20151009350).
J. Wolf has been supported by the German Research Foundation (DFG) through the project WO1988/1-1; 612414.  
M. Yang has been supported by the National Research Foundation of Korea (NRF) grant funded by the Korea government(MSIP) (No. 2016R1C1B2015731).

\appendix

\section{Proof of Lemma \ref{L22}}
\label{SA1}

The proof of Lemma \ref{L22} relies on the following proposition.

\begin{prop}
\label{prop10}
Let $ u\in V^2(Q_{1})$ be a local suitable weak solution to the Navier-Stokes equations.  We define 
$ v= u+ \nabla p_h $, where $ \nabla p_h = - E^{\ast}_{ B_{3/4}}(u)$. 
There exist absolute positive numbers $ K_{\ast}$ and $ \zeta $ such that if 
\begin{equation}
\int_{Q_{1}} | u|^3 dz  \le \zeta ^3
\label{A.11a}
\end{equation}
then for all $z_0 \in Q_{1/2}$ and for all natural number $k \ge 2$
\begin{equation}
\label{AG}
\fint_{Q(z_0, r_k)} |v|^3 dz \le K_{\ast}^3 \int_{Q_{1}} |u|^3 dz
\end{equation}
where $r_k = 2^{-k}$.
\end{prop}

We postpone the proof of Proposition \ref{prop10} at Appendix \ref{SA2}.
Suppose the proposition holds true.
Then using the Lebesgue differentiation theorem and \eqref{AG} we obtain that for almost all $z_0=(x_0, t_0) \in Q_{1/2}$
\begin{equation} 
|v(x_0,t_0)| \le K_{\ast} \left(\int_{Q_{1}} | u|^3 dz\right)^{1/3}.
\label{A.20}
\end{equation}
Using the triangular inequality and the mean value property of harmonic functions, we conclude that for almost all $ (x_0,t_0) \in Q(0,1/2)$
\begin{align*}
|u(x_0,t_0)| &\le |v(x_0,t_0)| + |\nabla p_h(x_0,t_0)| \\
&\le K_{\ast} \left(\int_{Q_{1}} | u|^3 dz\right)^{1/3} + 
c \| u(t_0)\|_{L^2(B_1)}. \\
&\le  K_{\ast} \left(\int_{Q_{1}} | u|^3 dz\right)^{1/3} + 
c \esssup_{t\in (-1,0)} \| u(t)\|_{ L^2(B_1)},
\end{align*} 
and hence
\begin{equation}
\| u\|_{ L^\infty(Q_{1/2})} \le K_{\ast} \left(\int_{Q_{1}} | u|^3 dz\right)^{1/3} + 
c \esssup_{t\in (-1,0)} \| u(t)\|_{ L^2(B_1)}. 
\label{A.21}
\end{equation}
Now, the assertion \eqref{E.28} in Lemma \ref{L22} follows from  \eqref{A.21} by a routine scaling argument.  
This completes the proof of Lemma \ref{L22}.

\section{Proof of Proposition \ref{prop10}}
\label{SA2}

We finally present the proof of Proposition \ref{prop10}.
The proof is divided into several steps.

\begin{itemize}
\item[\textbf{Step 1)}]
We shall prove the key inequality \eqref{AG} in Proposition \ref{prop10} by using a strong induction argument on $k$.
Let $K_{\ast}>1$ be a constant wihch will be specified at the final moment.
From the definition of a local suitable weak solution the following local energy inequality holds true for every nonnegative $\phi \in C^{\infty}_{\rm c} (B_{3/4} \times (-9/16,0])$ and almost 
all $s \in (-9/16,0]$ 
\begin{equation}
\begin{split}
\label{A.11b}
&\int |v(s)|^2 \phi(s) dx + 2 \int_{-r_3^2}^{s} \int |\nabla v|^2 \phi dz \\
&\le \int_{-r_3^2}^{s} \int |v|^2 (\partial_t + \Delta) \phi dz
+ \int_{-r_3^2}^{s} \int |v|^2 (v - \nabla p_h) \cdot \nabla \phi dz \\
&\quad + 2 \int_{-r_3^2}^{s} \int (v \otimes v - v \otimes \nabla p_h  : \nabla ^2 p_h) \phi dz \\
&\quad + 2 \int_{-r_3^2}^{s} \int (p_1+p_2) v \cdot \nabla \phi dz
\end{split}
\end{equation}
where 
\[\nabla p_1 = - E^{\ast}_{ B_{3/4}} (\divv (u \otimes u)), \qquad 
\nabla p_2 = E^{\ast}_{ B_{3/4}}(\Delta u).\]
Note that $v=u-\nabla p_h$ and so 
\begin{equation}
u \otimes u =v \otimes v - v \otimes \nabla p_h - \nabla p_h \otimes v+ \nabla p_h \otimes \nabla p_h
\label{A.11c}
\end{equation} 
almost everywhere in $Q_{3/4}$.

\item[\textbf{Step 2)}]
It is readily seen that \eqref{A.11a} holds for $k=2$.
Assume \eqref{AG} is true for $k = 2, \ldots, n$. 
Let $z_0\in Q_{1/4}$ be arbitrarily chosen and 
\[r_{ n+1} \le r \le  r_3.\]
Using the Cauchy--Schwarz inequality, the inductive assumption, and the fact that $p_h$ is harmonic, we get 
\begin{equation}
\begin{split}
\label{A.12}
\fint_{Q(z_0, r)} |v|^{3/2} |\nabla p_h|^{3/2} dz 
&\le \left(\fint_{Q(z_0, r)} |v|^3 dz\right)^{1/2}
\left(\fint_{Q(z_0, r)} |\nabla p_h|^3 dz\right)^{1/2} \\
&\lesssim r^{-5/2} K_{\ast}^{3/2}  \left(\int_{Q_{1}}  
|u|^{3} dz\right)^{1/2}\left(\int_{Q(z_0, r)}  | \nabla p_h|^{3} dz\right)^{1/2}
\cr
&\lesssim r^{ -1} K_{\ast}^{3/2} \int_{Q_{1}}  
|u|^{3} dz.
\end{split}
\end{equation}
Furthermore, applying the Poincar\'e inequality and using properties of harmonic functions, we find 
\begin{equation}
\begin{split}
&\fint_{Q(z_0, r)}  | \nabla p_h \otimes \nabla p_h - \inn{\nabla p_h \otimes \nabla p_h}_{B(x_0,r)}|^{3/2} dz \\
&\lesssim r^{-5 + 3/2}  \int_{Q(z_0, r)}  | \nabla p_h|^{3/2} | \nabla^2 p_h|^{3/2} dz \\
&\lesssim r^{-1/2}\int_{Q(0, 3/4)} | \nabla p_h|^{3} dz \\
&\lesssim r^{-1/2} \int_{Q_{1}} |u|^{3} dz. 
\label{A.13}
\end{split}
\end{equation}
Using the identity \eqref{A.11c} and combining the inductive assumption \eqref{AG} with the estimates \eqref{A.12} and \eqref{A.13}, we obtain that for all $ r_{ n+1} \le  r \le 1$
\[\int_{Q(z_0, r) }| u \otimes  u - \inn{u \otimes  u}_{B(x_0, r)} |^{3/2} dz
\lesssim K^3_{ \ast}r^{ 4}\int_{Q_{1}} |u|^{3} dz.  
\]
Applying Lemma\,2.8 in \cite{ChaeWolf2016a}, we find that for all $r_{ n+1} \le  r \le r_2$
\begin{equation}
\int_{Q(z_0, r)} | p_{ 1} - \inn{p_1}_{B(x_0, r)}|^{3/2} dz
\lesssim K_{\ast}^3  r^{ 4}\int_{Q_{1}}  |u|^{3} dz.  
\label{A.14}
\end{equation}

\item[\textbf{Step 3)}]
We denote by $\Psi_{n+1}$ the fundamental solution of the backward heat equation 
having its singularity at $ (x_0, t_0 + r_{ n+1}^2)$.
More precisely, for $(x, t) \in \R^{3}\times (-\infty, t_0+r_{ n+1})$ 
\[\Psi_{n+1} (x,t) = \frac{c_0}{(r_{ n+1}^2- t+ t_0 )^{3/2}} 
\exp \Set{- \frac{| x-x_0|^2}{4(r^2_{ n+1}- t+ t_0)}}.\]
Taking a suitable cut off function $ \chi \in C^{\infty}( \R^{n})$ for $Q(z_0, r_4)\subset  Q(z_0, r_3)$, we may insert 
$\Phi_{n+1} := \Psi_{n+1} \chi $ into the local energy inequality \eqref{A.11b}  to get for almost all $s \in (t_0-r_3^2, t_0)$
\begin{align*}
&\int_{B(x_0, r_3)} \Phi_{n+1}(s) |v(s)|^2 dx 
+ 2 \int_{t_0-r_3^2}^{s}\int_{B(x_0, r_3)}\Phi_{n+1} | \nabla v|^2 dz \\ 
&\le \int_{t_0-r_3^2}^{s} \int_{B(x_0, r_3)} |v|^2 (\partial_t + \Delta) \Phi_{n+1} dz \\
&\quad + \int_{t_0-r_3^2}^{s}\int_{B(x_0, r_3)} |v|^2 (v-\nabla p_h) \cdot \nabla \Phi_{n+1} dz \\
&\quad + 2 \int_{t_0-r_3^2}^{s}\int_{B(x_0, r_3)} (v \otimes v - v \otimes \nabla p_h : \nabla ^2 p_h) \Phi_{n+1} dz \\
&\quad  + 2 \int_{t_0-r_3^2}^{s}\int_{B(x_0, r_3)} (p_1+p_2) v \cdot \nabla \Phi_{n+1} dz.
\end{align*}
Arguing as in \cite{MR673830}, we obtain from the above inequality that 
\begin{equation}
\begin{split}
&\esssup_{s \in (t_0-r_{ n+1}^2, t_0)} \fint_{B(x_0, r_{ n+1})} |v(s)|^2 dx 
+ r_{ n+1}^{ -3}\int_{Q(z_0, r_{ n+1})}  | \nabla v|^2 dz \\ 
&\lesssim \int_{Q(z_0, r_3)} |v|^2 |(\partial_t + \Delta) \Phi_{n+1}| dz
+ \int_{Q(z_0, r_3)} |v|^2 (|v|+|\nabla p_h|) |\nabla \Phi_{n+1}| dz \\ 
&\quad + \int_{Q(z_0, r_3)} |v| (|v|+|\nabla p_h|) |\nabla^2 p_h| \Phi_{n+1} dz 
+ \int_{Q(z_0, r_3)}  (p_1+p_2)  v\cdot  \nabla \Phi_{n+1} dz \\
&=: I_1 + I_2 + I_3 + I_4.
\label{A.16}
\end{split}
\end{equation}

\item[\textbf{Step 4)}]
In this step we shall estimate the integrals $I_1$, $I_2$, and $I_3$.
They can be handled by the similar way.

Obviously, we have $|(\partial_t + \Delta) \Phi_{n+1}| \le C$ in $ Q(z_0, r_3)$ so that 
\[I_1 \le C \| v\|^2_{ L^3(Q(z_0, r_3)} \lesssim \left(\int_{Q_{1}}  |u|^{3} dz\right)^{2/3}.\]

Using $|\nabla \Phi_{n+1}| \le Cr^{ -4}_k$ in $ Q(z_0, r_k) \setminus Q(z_0, r_{ k+1})$ for all $ k=2, \dots, n$ and the inductive assumption \eqref{AG}, we obtain 
\begin{align*}
&\int_{Q(z_0, r_3)} |v|^3 |\nabla \Phi_{n+1}| dz \\
&= \sum_{k=3}^{n} \int_{Q(z_0, r_{ k}) \setminus Q(z_0, r_{ k+1})} |v|^3  | \nabla \Phi_{n+1}| dz
+ \int_{Q_{ r_{ n+1}}(z_0) } |v|^3  | \nabla \Phi_{n+1}| dz \\
&\lesssim K_{\ast}^3 \sum_{k=2}^{n} r_k^{ -4} r_k^{ 5} \int_{Q_{1}}  |u|^{3} dz
\lesssim K_{\ast}^3\int_{Q_{1}}  |u|^{3} dz. 
\end{align*}

Similarly,
\begin{align*}
&\int_{Q(z_0, r_3)} |v|^2 |\nabla p_h| |\nabla \Phi_{n+1}| dz \\
&= \sum_{k=3}^{n} \int_{Q(z_0, r_{ k})  \setminus Q(z_0, r_{ k+1})} |v|^2 |\nabla p_h|  | \nabla \Phi_{n+1}| + 
\int_{Q(z_0, r_{ n+1})  } |v|^2 | \nabla p_h|  | \nabla \Phi_{n+1}| 
\\
&\lesssim K_{\ast}^2 \sum_{k=1}^{n} r_k^{ -4}  r_k^{13/3} \int_{Q_{1}}  |u|^{3} dz 
\lesssim K_{\ast}^3 \int_{Q_{1}}  |u|^{3} dz. 
\end{align*}
Hence we have $I_2 \lesssim K_{\ast}^3 \int_{Q_{1}}  |u|^{3} dz$ and the implied constant does not depend on $n$.

Using $\Phi_{n+1} \le Cr^{ -3}_k$ in $ Q(z_0, r_{ k}) \setminus Q(z_0, r_{ k+1})$ for all $ k=1, \ldots, n+1$, the inductive assumption \eqref{AG}, and the properties of harmonic functions, we get 
\begin{align*}
&\int_{Q(z_0, r_3)} |v|^2 |\nabla^2 p_h| \Phi_{n+1} dz \\
&= \sum_{k=3}^{n} \int_{Q(z_0, r_{ k})  \setminus  Q(z_0, r_{ k+1})} |v|^2 | \nabla^2 p_h| \Phi_{n+1} dz
+ \int_{Q(z_0, r_{ n+1}) } |v|^2 | \nabla^2 p_h|  | \Phi_{n+1}| dz \\
&\lesssim K_{\ast}^2 \sum_{k=2}^{n} r_k^{ -3}  r_k^{13/3} \int_{Q_{1}}  |u|^{3} dz
\lesssim K_{\ast}^3 \int_{Q_{1}}  |u|^{3} dz. 
\end{align*}
Similarly,
\begin{align*}
&\int_{Q(z_0, r_3)} |v| |\nabla p_h| |\nabla^2 p_h| \Phi_{n+1} dz \\
&= \sum_{k=3}^{n} \int_{Q(z_0, r_{ k})  \setminus Q(z_0, r_{ k+1})} |v| |\nabla p_h|| \nabla^2 p_h |  \Phi_{n+1} dz 
+ \int_{Q(z_0, r_{ n+1})} |v| |\nabla p_h|| \nabla^2 p_h |  \Phi_{n+1} dz \\
&\lesssim K_{\ast} \sum_{k=2}^{n} r_k^{ -3}  r_k^{11/3} \int_{Q_{1}}  |u|^{3} dz
\lesssim K_{\ast}^3 \int_{Q_{1}}  |u|^{3} dz. 
\end{align*}
Hence we have $I_3 \lesssim K_{\ast}^3 \int_{Q_{1}}  |u|^{3} dz$ and the implied constant does not depend on $n$.

\item[\textbf{Step 5)}]
In this step we estimate the last integral $I_4$ in \eqref{A.16}.
We argue as in \cite{MR673830}. 
Let $ \chi _k$ denote cut-off functions, suitable for $Q(z_0, r_{ k+1}) \subset Q(z_0, r_{ k})$, $k=3, \ldots, n+1$. 
Since $v$ is divergence free, we can subtract an average from $p_2$ and use the partition of unity so that 
\begin{align*}
&\int_{Q(z_0, r_3)} p_2 v\cdot  \nabla \Phi_{n+1} dz \\
&= \sum_{k=3}^{n} \int_{Q(z_0, r_{ k})  \setminus Q(z_0, r_{ k+2})} (p_2 - \inn{p_2}_{B(x_0,r_k)}) v\cdot \nabla (\Phi_{n+1} (\chi _k - \chi _{ k+1})) \\
&\quad + \int_{Q(z_0, r_{2})} p_2  v\cdot \nabla (\Phi_{n+1} (1- \chi _3)) \\
&\quad + \int_{Q(z_0, r_{ n+1})} (p_2 - \inn{p_2}_{B(x_0,r_{n+1})}) v\cdot \nabla (\Phi_{n+1} \chi _{ n+1}) \\
&=: J_1 + J_2 + J_3.  
\end{align*}
As $ | \nabla (\Phi_{n+1}(\chi_k - \chi_{ k+1}))|  \le C  r_k^{ -4}$ for $ k=1, \ldots, n$, applying Poincar\'e's inequality, 
using the fact that  $ p_2 $ is harmonic, together  with (57)$ _k$ and \eqref{E26}  we see that 
\begin{align*}
&\int_{Q(z_0, r_{ k})  \setminus Q(z_0, r_{ k+2})} (p_2 - \inn{p_2}_{B(x_0,r_{k})}) v\cdot \nabla (\Phi_{n+1} (\chi _k - \chi _{ k+1})) \\
&\lesssim K_{\ast}  r_k^{ -4} r_k^{ 5}  \left(\int_{Q_{1}} | u|^3 dz\right)^{1/3} \left(\int_{Q_{1/2}} p_2^2 dz\right)^{1/2}   
\\
&\lesssim K_{\ast}  r_k \left(\int_{Q_{1}} | u|^3 dz\right)^{1/3} 
\left(\int_{Q_{3/4}} | \nabla u|^2 dz\right)^{1/2} \\
&\lesssim K_{\ast} r_k \left(\int_{Q_{1}} | u|^3 dz\right)^{2/3}.   
\end{align*}
Summation from $ k=3 $ to $ n$ yields 
\[J_1 \lesssim K_{\ast} \left(\int_{Q_{1}} | u|^3 dz\right)^{2/3}.\]
Similarly, we can make 
\[J_2 + J_3 \lesssim (1+r_{ n+1}) K_{\ast} \left(\int_{Q_{1}} | u|^3 dz\right)^{2/3}.\]
Thus,
\[\int_{Q(z_0, r_3)} p_2 v\cdot  \nabla \Phi_{n+1} dz \lesssim K_{\ast} \left(\int_{Q_{1}} | u|^3 dz\right)^{2/3}.\]

Finally, arguing in the same way and making use of \eqref{A.14}, we can get 
\[\int_{Q(z_0, r_3)} p_1 v\cdot  \nabla \Phi_{n+1} dz \lesssim K_{\ast}^3 \int_{Q_{1}} | u|^3 dz\] 
and therefore 
\[I_4 \lesssim K_{\ast}^3 \int_{Q_{1}} | u|^3 dz + K_{\ast} \left(\int_{Q_{1}} | u|^3 dz\right)^{2/3}.\]

\item[\textbf{Step 6)}]
Inserting the estimates of $I_1$, $I_2$, $I_3$, and $I_4$ into the right-hand side of \eqref{A.16}, we obtain that for some absolute constant $C_1>0$, independently of $n$, 
\begin{align*}
&\esssup_{s \in (t_0-r_{ n+1}^2, t_0)} \fint_{B_{ r_{ n+1}}(x_0)} |v(s)|^2  
+ r_{ n+1}^{ -3} \int_{Q(z_0, r_{ n+1})}  | \nabla v|^2 dz \\ 
&\le C_1 \left(K_{\ast}^3 \int_{Q_{1}} | u|^3 dz + K_{\ast} \left(\int_{Q_{1}} |u|^3 dz\right)^{2/3}\right) \\
&= \left(C_1 K_{\ast} E + \frac{C_1}{K_{\ast}}\right) K_{\ast}^2 E^2
\end{align*}
where $E=\left(\int_{Q_{1}} | u|^3 dz\right)^{1/3}$. 

On the other hand, using a standard interpolation, we obtain that 
\[\fint_{Q(z_0, r_{ n+1})} |v|^3 dz
\le C_2 \left(\esssup_{s \in (t_0-r_{ n+1}^2, t_0)}
\fint_{B(x_0, r_{ n+1})} |v(s)|^2 + r_{ n+1}^{ -3}\int_{Q(z_0, r_{ n+1})}  | \nabla v|^2\right)^{3/2}\]
for some absolute constant $C_2>1$ and hence
\[\fint_{Q(z_0, r_{ n+1})} |v|^3 dz 
\le \left(C_2C_1 K_{\ast} E + \frac{C_2C_1}{K_{\ast}}\right)^{3/2} K_{\ast}^3
E^3.\]
Note that neither $ C_1$ nor  $C_2 $ depend on the choice of $ K_{\ast}$.  
Thus, we may set 
\[K_{\ast} = 2C_1C_2,\qquad  \zeta  = \frac{1}{4C_1^2C_2^2}\]
so that if $ E \le \zeta $, then
\[\fint_{Q(z_0, r_{ n+1})} |v|^3 \le K_{\ast}^3 E^3 = K_{\ast}^3 \int_{Q_{1}} | u|^3 dz.\]
Hence \eqref{AG} is true for $k=n+1$.
\end{itemize}
This completes the proof of Proposition \ref{prop10}.


\begin{thebibliography}{30}

\bibitem{MR673830}
L.~Caffarelli, R.~Kohn, and L.~Nirenberg.
\newblock Partial regularity of suitable weak solutions of the
  {N}avier-{S}tokes equations.
\newblock {\em Comm. Pure Appl. Math.}, 35(6):771--831, 1982.

\bibitem{ChaeWolf2016a}
D.~Chae and J.~Wolf.
\newblock On the liouville theorem for self similar solutions to the {N}avier-{S}tokes equations.
\newblock {\em Preprint}, 2016.

\bibitem{MR1780481}
H.~J. Choe and J.~L. Lewis.
\newblock On the singular set in the {N}avier-{S}tokes equations.
\newblock {\em J. Funct. Anal.}, 175(2):348--369, 2000.

\bibitem{MR2005639}
L.~Escauriaza, G.~Seregin, and V.~{\v{S}}ver{\'a}k.
\newblock Backward uniqueness for parabolic equations.
\newblock {\em Arch. Ration. Mech. Anal.}, 169(2):147--157, 2003.

 \bibitem{gal}G. Galdi, C. Simader and H. Sohr, {\it On the Stokes problem in Lipshitz domain,}  Annali di Mat. Pura Appl., (IV), 167 (1994), pp. 147-163.

\bibitem{MR2449250}
L. Grafakos.
\newblock {\em Classical and modern {F}ourier analysis}.
\newblock Pearson Education, Inc., Upper Saddle River, NJ, 2004.

\bibitem{MR2308753}
S. Gustafson, K. Kang, and T. P. Tsai.
\newblock Interior regularity criteria for suitable weak solutions of the
  {N}avier-{S}tokes equations.
\newblock {\em Comm. Math. Phys.}, 273(1):161--176, 2007.

\bibitem{MR0050423}
E. Hopf.
\newblock \"{U}ber die {A}nfangswertaufgabe f\"ur die hydrodynamischen
  {G}rundgleichungen.
\newblock {\em Math. Nachr.}, 4:213--231, 1951.

\bibitem{MR1632780}
H. Kozono.
\newblock Removable singularities of weak solutions to the {N}avier-{S}tokes equations.
\newblock {\em Comm. Partial Differential Equations}, 23:949--966, 1998.

\bibitem{MR0236541}
O.~A. Lady{\v{z}}enskaja.
\newblock Uniqueness and smoothness of generalized solutions of
  {N}avier-{S}tokes equations.
\newblock {\em Zap. Nau\v cn. Sem. Leningrad. Otdel. Mat. Inst. Steklov.
  (LOMI)}, 5:169--185, 1967.

\bibitem{MR1738171}
O.~A. Ladyzhenskaya and G.~A. Seregin.
\newblock On partial regularity of suitable weak solutions to the
  three-dimensional {N}avier-{S}tokes equations.
\newblock {\em J. Math. Fluid Mech.}, 1(4):356--387, 1999.

\bibitem{MR1555394}
J. Leray.
\newblock Sur le mouvement d'un liquide visqueux emplissant l'espace.
\newblock {\em Acta Math.}, 63(1):193--248, 1934.

\bibitem{MR1488514}
F. Lin.
\newblock A new proof of the {C}affarelli-{K}ohn-{N}irenberg theorem.
\newblock {\em Comm. Pure Appl. Math.}, 51(3):241--257, 1998.

\bibitem{MR1738173}
J. Neustupa.
\newblock Partial regularity of weak solutions to the {N}avier-{S}tokes
  equations in the class {$L^\infty(0,T;L^3(\Omega)^3)$}.
\newblock {\em J. Math. Fluid Mech.}, 1(4):309--325, 1999.

\bibitem{MR0126088}
G. Prodi.
\newblock Un teorema di unicit\`a per le equazioni di {N}avier-{S}tokes.
\newblock {\em Ann. Mat. Pura Appl. (4)}, 48:173--182, 1959.

\bibitem{MR0454426}
V. Scheffer.
\newblock Partial regularity of solutions to the {N}avier-{S}tokes equations.
\newblock {\em Pacific J. Math.}, 66(2):535--552, 1976.

\bibitem{MR0510154}
V. Scheffer.
\newblock Hausdorff measure and the {N}avier-{S}tokes equations.
\newblock {\em Comm. Math. Phys.}, 55(2):97--112, 1977.

\bibitem{MR1829531}
G. Seregin.
\newblock On the number of singular points of weak solutions to the
              {N}avier-{S}tokes equations.
\newblock {\em Comm. Pure Appl. Math.}, 8:1019--1028, 2001.

\bibitem{MR0150444}
J. Serrin.
\newblock The initial value problem for the {N}avier-{S}tokes equations.
\newblock In {\em Nonlinear {P}roblems ({P}roc. {S}ympos., {M}adison, {W}is.,
  1962)}, pages 69--98. Univ. of Wisconsin Press, Madison, Wis., 1963.

\bibitem{MR1877269}
H. Sohr.
\newblock A regularity class for the {N}avier-{S}tokes equations in
              {L}orentz spaces.
\newblock {\em J. Evol. Equ.}, 1(4):441--467, 2001.

\bibitem{WZ2014}
W.~Wang and Z.~Zhang.
\newblock On the interior regularity criteria and the number of singular points to the Navier--Stokes equations.
\newblock {\em J. Anal. Math.}, 123(1):139--170, 2014. 

\bibitem{Wolf2015c}
J.~Wolf. 
\newblock On the local regularity of suitable weak solutions to the generalized {N}avier-{S}tokes equations.
\newblock {\em Ann. Univ. Ferrara}, 61:149--171, 2015.

\bibitem{Wolf2016}
J.~Wolf. 
\newblock On the local pressure of the {N}avier-{S}tokes equations and related system.
\newblock {\em to appear in Adv. Diff. Equs.}, 2016.

\end{thebibliography}
\end{document}